\documentclass[10pt,english]{smfart}
\usepackage{graphics,amssymb,amsfonts,amsmath}
\usepackage[all,ps,cmtip]{xy} 
\usepackage{smfthm,bull} 

\usepackage{mathrsfs} 
\let\mathcal\mathscr

\def\ie{i.e.}
\def\eg{e.g.}
\def\isom{\simeq}

\def\codim{\mathop{\rm codim}\nolimits}

\def\Ker{\mathop{\rm Ker}\nolimits}

\def\Supp{\mathop{\rm Supp}\nolimits}

\def\Id{\mathop{\rm Id}\nolimits}

\def\Sym{\mathbf{Sym}}

\def\iff{if and only if}

\def\lra{\longrightarrow}
\def\llra{\hbox to 12mm{\rightarrowfill}}

\def\vide{\varnothing}

\def\av{abelian variety}

\def\Pic{\mathop{\rm Pic}\nolimits}

\def\Proj{\mathop{\rm Proj}\nolimits}
\def\tr{\mathop{\rm Tr}\nolimits}

\def\PA{\Pic^0(A)}

\def\PC{\Pic^0(C)}

\def\PX{\Pic^0(X)}
\def\PY{\Pic^0(Y)}
\def\Pa{\Pic^0(\alpha)}
\def\Pf{\Pic^0(f)}
\def\Pg{\Pic^0(g)}
\def\Pp{\Pic^0(\pi)}

\def\Z{{\bf Z}}

\def\P{{\bf P}}

\def\C{{\bf C}}

\def\L{\Lambda}

\def\cF{{\mathcal F}}
\def\cG{{\mathcal G}}

\def\cO{{\mathcal O}}

\def\cS{{\mathcal S}}

\def\cgg{continuously globally generated}

\def\gge{globally generated}

\theoremstyle{plain}
\author{Olivier Debarre}
\email{debarre@math.u-strasbg.fr}
\urladdr{http://www-irma.u-strasbg.fr/~debarre/}
\title{On coverings of simple abelian varieties}
\alttitle{Sur les rev\^etements des vari\'et\'es ab\'eliennes simples}

\begin{document} 
\frontmatter
\begin{abstract}To any finite covering $f:Y\to X$ of degree $d$ between smooth complex projective manifolds, one associates a vector bundle $E_f$ of rank $d-1$ on $X$ whose total space contains $Y$. 
It is known that $E_f$ is ample when $X$ is a projective space (\cite{la}), a Grassmannian (\cite{man}), or a lagrangian Grassmannian (\cite{km}). %
 We show an analogous result when $X$ is a simple abelian variety and $f$ does not factor through any nontrivial isogeny $X'\to X$. This result is obtained by showing that $E_f$ is $M$-regular in the sense of Pareschi--Popa, and that any $M$-regular sheaf is ample.
\end{abstract}

\begin{altabstract}
\`A tout rev\^etement fini $f:Y\to X$ de degr\'e $d$ entre vari\'et\'es projectives lisses complexes, on associe un fibr\'e vectoriel $E_f$ de rang $d-1$ sur $X$ dont l'espace total contient $Y$. 
On sait que $E_f$ est ample lorsque $X$ est un espace projectif (\cite{la}),   une grassmannienne (\cite{man}) ou une grassmannienne lagrangienne (\cite{km}). %
 Nous montrons un r\'esultat analogue lorsque $X$ est une vari\'et\'e ab\'elienne simple et que $f$ ne se factorise par aucune isog\'enie non triviale $X'\to X$. Ce r\'esultat est obtenu en montrant que $E_f$ est $M$-r\'egulier au sens de Pareschi--Popa, puis que tout faisceau $M$-r\'egulier est ample.
\end{altabstract}

\altkeywords{Vari\'et\'e ab'elienne, fibr\'e vectoriel, faisceau ample, faisceau $M$-r\'egulier, faisceau contin\^ument engendr\'e, th\'eor\`eme de Barth--Lefschetz, transform\'ee de Mukai}

 \maketitle  

\mainmatter
\section{Introduction}

We work over the complex numbers. Let $f:Y\to X$ be a finite surjective morphism of degree $d$ between smooth projective varieties of the same dimension $n$.  The morphism $f$ is flat, hence the sheaf $f_*\cO_Y$ is locally free. We may define a locally free sheaf $E_f$ of rank $d-1$ on $X$ as the dual of the kernel of the trace map $\tr_{Y/X}:f_*\cO_Y\to\cO_X$, so that
$$f_*\cO_Y=\cO_X\oplus E_f^*$$
By duality for a finite flat morphism, we have
$$f_*\omega_{Y/X}=\cO_X\oplus E_f$$
Our aim  is to prove the following statement conjectured in \cite{de}.

\begin{theo}
Let $X$ be a {\em simple} abelian variety, let $Y$ be a smooth {\em connected} projective variety,
and let $f:Y\to X$ be a finite cover. If $f$  does not factor through any nontrivial isogeny $X'\to X$, the vector bundle  $E_f$ is ample.
\end{theo}

For a more general statement, see Theorem \ref{main}. See also the remarks at the end of this article for more
comments.  Even if $X$ is not simple, the vector bundle $E_f$ is known   to be nef   (\cite{ps},
  Theorem 1.17; \cite{laz}, Example 6.3.59) and its restriction to a general complete intersection curve in $X$  to be ample (\cite{hkp}, Lemma 2.7).

  The ampleness of $E_f$ has a number
of consequences, as explained in \cite{laz}, Example 6.3.56. In our case,   one new
 statement beyond the Fulton--Hansen-type results already obtained in \cite{de} is the following:
under the hypotheses of the theorem,
 the induced morphism
$$H^i(f,\C):H^i(X,\C)\to H^i(Y,\C)
$$ is bijective for $i\le n-d+1$ (\cite{laz}, Theorem 7.1.16).

When moreover $d\le n$,   the morphism $ \pi_1(f):\pi_1(Y)\to  \pi_1(X)
$ is bijective.\footnote{For algebraic fundamental groups, this is  \cite{de},
 Corollaire 6.2; for topological fundamental groups, this is  \cite{livre}, Exercice VIII.5,
 where the hypothesis $d\le n$ is unfortunately missing.} In particular, the
 group $H_1(Y,\Z)$ is isomorphic to $ H_1(X,\Z)$, hence is torsion-free, and so
is  $H^2(Y,\Z)$ by the universal coefficient theorem. 

When $d\le n-1$, the morphism $H^2(f,\Z):H^2(X,\Z)\to H^2(Y,\Z)
$ is   injective with finite cokernel, hence so is
$\Pic(f):\Pic(X)\to\Pic(Y)$. %
  It seems likely that those two maps are bijective.

The proof is a simple application of the results of \cite{pp1} about global generation of sheaves on an abelian variety. More precisely, it is based on the remark that any $M$-regular sheaf (\S~\ref{s2}) on an abelian variety is ample (Corollary \ref{coro}).

This work was done while the author was visiting the University of Michigan at Ann Arbor. Many thanks to Bill Fulton and Rob Lazarsfeld for support and many stimulating conversations.

\section{Ample sheaves}\label{amp}

To any coherent sheaf $\cF$ on a scheme $X$ of finite type  over $\C$, one associates the $X$-scheme
$$\P(\cF)=\Proj\Big(\bigoplus_{m\ge 0}\Sym^m\cF\Big)
$$
and an invertible sheaf $\cO_{\P(\cF)}(1)$ on $\P(\cF)$. The sheaf $\cF$ is said to be ample if $\cO_{\P(\cF)}(1)$ is.

 Well-known properties of ampleness for locally free sheaves (see for example \cite{laz}, Chapter 6) still hold in this general setting:

a) the sheaf $\cF$ is ample \iff, for any coherent sheaf $\cG$ on $X$, the sheaf $\cG\otimes\Sym^m\cF$ is globally generated for all $m\gg 0$ (\cite{ku}, Theorem 1);

b) any   quotient of an ample sheaf is ample (\cite{ku}, Proposition 1);

c) if $\pi:Y\to X$ is a finite morphism,  $\cF$ is ample \iff\ $\pi^*\cF$ is (this is because $\P(\pi^*\cF)=\P(\cF)\times_XY$ and $\cO_{\P(\cF)}(1)$ pulls back, by a finite morphism, to  $\cO_{\P(\pi^*\cF)}(1)$);

d) if $X$ is proper and $\cF$ is globally generated, $\cF$ is ample \iff, for any curve $C$ in $X$, the restriction $\cF\otimes\cO_C$ has no trivial quotient (Gieseker's Lemma).

\section{Continuously generated sheaves}\label{s2}

Following \cite{pp1}, Definition 2.10, we say that a coherent sheaf $\cF$ on
an  irreducible projective variety $X$ is {\em \cgg} if, for any nonempty subset $U$ of $\PX$, the sum of the twisted evaluation maps
$$\bigoplus_{\xi\in U} H^0(X,\cF\otimes P_\xi)\otimes P^\vee_\xi \to \cF$$
is surjective, where, for any element $\xi$ of $\Pic^0(X)$, we denote by $P_\xi$ the corresponding numerically trivial line bundle on $X$. This property is equivalent to the existence of a positive integer $N$ such that for
 $(\xi_1,\dots,\xi_N)$   general in $\PX^N$, the analogous map
\begin{equation}\label{def}
\bigoplus_{i=1}^N H^0(X,\cF\otimes P_{\xi_i})\otimes P_{\xi_i}^\vee\to  \cF 
\end{equation}
is surjective. Being a quotient of a direct sum of numerically trivial line bundles, a  \cgg\ sheaf  is nef. Our aim is to show that under certain circumstances, it is ample.

\begin{prop}\label{lemm}
A  coherent sheaf $\cF$  on an irreducible projective variety  $X$ is \cgg\ \iff\
there exists a connected abelian
Galois \'etale cover $\pi:Y\to X$ such that $\pi^*(\cF\otimes P_\xi)$ is \gge\ for all $\xi\in\PX$.
\end{prop}

\begin{proof} Assume $\cF$ is \cgg\ and let ${\xi_0}\in\PX$. Since torsion points are dense in $\PX^N$, the open subset of $\PX^N$ of points for which the map (\ref{def}) is surjective and all $h^0(X,\cF\otimes P_{\xi_i})$ are minimal contains a point of the type 
$$({\xi_0}+\eta_1({\xi_0}),\dots,{\xi_0}+\eta_N({\xi_0}))$$
where $(\eta_1({\xi_0}),\dots,\eta_N({\xi_0}))$ is torsion, hence contains also $U_{\xi_0}+(\eta_1({\xi_0}),\dots,\eta_N({\xi_0}))$, where $U_{\xi_0}$ is a neighborhood of $\xi_0$ in $\PX$. Since $\PX$ is quasi-compact, it is covered by finitely many such neighborhoods, say $U_{\xi_1},\dots,U_{\xi_M}$.

Let 
$\pi:Y\to X$ be a connected abelian Galois \'etale cover such that the kernel of    $\Pp:\PX\to \PY$ contains all $\eta_i(\xi_j)$, for $i\in\{1,\dots,N\}$ and $j\in\{1,\dots,M\}$. Fix $j \in\{1,\dots,M\}$; the map
$$\bigoplus_{i=1}^N H^0(X,\cF\otimes P_\xi\otimes P_{\eta_i(\xi_j)})\otimes \pi^*P_\xi^\vee\otimes\pi^*P^\vee_{\eta_i(\xi_j)}\lra\pi^*\cF$$
is surjective for all $\xi\in U_{\xi_j}$. But this map is
$$\bigoplus_{i=1}^N H^0(X,\cF\otimes P_\xi\otimes P_{\eta_i(\xi_j)})\otimes \pi^*P_\xi^\vee\lra \pi^*\cF$$
and since each $H^0(X, \cF\otimes P_\xi\otimes P_{\eta_i(\xi_j)})$ is a vector subspace of $H^0(Y,\pi^*(\cF\otimes P_\xi))$,
the sheaf $\pi^*(\cF\otimes P_\xi)$ is globally generated for all $\xi\in U_{\xi_j}$, hence for all $\xi\in\PX$.

For the converse, assume that there exists a connected abelian Galois \'etale cover $\pi:Y\to X$ such that  the evaluation map
$$H^0(Y,\pi^*(\cF\otimes P_\xi))\otimes \cO_Y\to \pi^*(\cF\otimes P_\xi)
$$
is surjective for all $\xi\in\PX$. 
Since $\pi$ is finite, the map
$$ H^0(X,\cF\otimes P_\xi\otimes \pi_*\cO_Y) \otimes \pi_*\cO_Y
\to  \cF\otimes P_\xi\otimes \pi_*\cO_Y
$$
is also surjective. If we let $\Ker(\Pp)=\{\eta_1,\dots,\eta_N\}$, we have $\pi_*\cO_Y=\bigoplus_{i=1}^NP_{\eta_i}$,   the map
$$
\begin{matrix}
\left( \bigoplus_{i=1}^N H^0(X,\cF\otimes P_\xi\otimes P_{\eta_i})\right)  \otimes \left( \bigoplus_{i=1}^N   P_{\eta_i}\right)\\
\downarrow\\  \cF\otimes P_\xi\otimes \left( \bigoplus_{i=1}^N   P_{\eta_i}\right)
\end{matrix}
$$
is   surjective, and so is
$$
 \bigoplus_{i=1}^N H^0(X,\cF\otimes P_\xi\otimes P_{\eta_i})  \otimes   P_{\eta_i}^\vee\to  \cF\otimes P_\xi$$
In other words,   the map (\ref{def}) is surjective for  $(\xi_1,\dots,\xi_N)=(\xi+\eta_1,\dots,\xi+\eta_N)$, for all $\xi\in\PX$. Choosing $\xi_0$ such that $h^0(X,\cF\otimes P_{\xi_0+ \eta_i})$ takes the general (minimal) value for each $i$ in $\{1,\dots,N\}$, we obtain that the map (\ref{def}) is still surjective for 
$(\xi_1,\dots,\xi_N)$ in a neighborhood of $(\xi_0+\eta_1,\dots,\xi_0+\eta_N)$. This proves that $\cF$ is \cgg.\end{proof}

 \begin{coro}\label{coro}
Let $X$ an irreducible projective variety with a finite map to an abelian variety. Any
   \cgg\ coherent sheaf  on $X$ is ample.
\end{coro}

The converse is in general false: if $L$ is an ample line bundle on an abelian variety $A$ of dimension $g$, a general map $(L^{-d})^{\oplus g}\to (L^{-1})^{\oplus 2g}$ is injective for $d\gg0$ and its cokernel is an ample vector bundle $E$ (\cite{laz}, Theorem 6.3.65). If $g\ge 2$, we have $H^0(A,E\otimes P_\xi)=0$ for all $\xi\in\PA$, hence $E$ cannot be \cgg.

\begin{proof} Let $\cF$ be a \cgg\ coherent sheaf on $X$. By Proposition \ref{lemm}, there exists a connected abelian Galois \'etale cover $\pi:Y\to X$ such that $\pi^*(\cF\otimes P_\xi)$ is \gge\ for all $\xi\in\PX$. 

Let $C$ be a curve in $Y$. If there is a trivial quotient $\pi^*\cF\vert_C\twoheadrightarrow\cO_C$, we have also surjections $\pi^*(\cF\otimes P_\xi)\vert_C\twoheadrightarrow \pi^*  P_\xi \vert_C$ for each $\xi\in\PX$. Since $\pi^*(\cF\otimes P_\xi)$ is \gge, so is $\pi^*  P_\xi \vert_C$. This implies that the composition $\PX\to\PY\to\PC$ is zero, hence that $\pi(C)$ is contracted by any map from $X$ to an abelian variety. This contradicts our hypothesis, hence 
$\pi^*\cF\vert_C$ has no trivial quotient.

By Gieseker's Lemma, $\pi^*\cF$ is ample, and so is $\cF$ (\S\ \ref{amp}).
\end{proof}

\section{The main theorem}

Following \cite{pp1}, Definition 2.1, we say that a coherent sheaf $\cF$ on an
 abelian variety $A$ is {\em $M$-regular} if
$$\codim_{\PA}\Supp \bigl(R^i\hat\cS(\cF)\bigr)>i
$$
for all $ i>0$ ($R^i\hat\cS$ is the $i$th Fourier--Mukai functor). This is the
case if
$$\codim_{\PA}\{\xi\in\PA\mid H^i(A,\cF\otimes P_\xi)\ne0 \}>i
$$
for all $i>0$. We
refer to \cite{muk} and
\cite{pp1} for more details. For our purposes, the main result of \cite{pp1} (Proposition
2.13)  is that {\em an $M$-regular coherent sheaf on an abelian variety is \cgg.}

\begin{theo}\label{main}
Let $X$ be a smooth connected projective variety with a finite map to a {\em
simple}
abelian variety, let $Y$ be a smooth {\em connected} projective variety with a finite surjective map  $f:Y\to X$. If $f$  factors through no nontrivial connected abelian Galois \'etale covering of $X$, the vector bundle  $E_f\otimes\omega_X$ is  ample.
\end{theo}

\begin{proof}Let $n$ be the common dimension of $X$ and $Y$, and let
$\alpha:X\to A$ be a finite map to a simple abelian variety  such that
$\Pa:\PA\to\PX$ is injective. Set $g=\alpha\circ f$. By \cite{gl1},
Theorem 1, \cite{gl2}, Theorem 0.1, and \cite{el}, Remark 1.6 (see also
\cite{el}, Theorem 1.2),  every irreducible component of
 the set
$$V_i=\{\xi\in\PA\mid H^{n-i}(Y,g^*P_\xi^\vee)\ne0\}$$ is a
translated abelian subvariety of $\PA$ of codimension at least $i$. In particular, since $A$ is simple,  $V_i$ is finite for $i>0$.

Since $Y$ is connected, we have
\begin{eqnarray*}
V_n&=&\{\xi\in\PA\mid H^0(Y,g^*P_\xi^\vee)\ne0\}\\
&=&\{\xi\in\PA\mid g^*P_\xi^\vee\isom\cO_Y\}\\
&=&\Ker(\Pg:\PA\to \PY\}
\end{eqnarray*}
hence $V_n= \{0\}$ since both $\Pa$ and $\Pf$ are injective ($f$  factors through no nontrivial abelian \'etale covering of $X$). Consider now
\begin{eqnarray*}
W_i&=&\{\xi\in\PA\mid H^i(X, E_f\otimes\omega_X \otimes \alpha^* P_\xi)\ne 0\}\\
&=&\{\xi\in\PA\mid H^i(A, \alpha_*(E_f\otimes\omega_X) \otimes  P_\xi)\ne 0\}
\end{eqnarray*}
By Serre duality on $Y$, 
\begin{eqnarray*}
V_i
&=&\{\xi\in\PA\mid H^i(Y,\omega_Y\otimes g^*P_\xi)\ne0\}\\
&=&\{\xi\in\PA\mid H^i(X,f_*\omega_Y\otimes \alpha^*P_\xi)\ne0\}
\end{eqnarray*}
Since
$f_*\omega_Y=f_*\omega_{Y/X}\otimes \omega_X=\omega_X\oplus ( E_f\otimes \omega_X)$,
we have $W_i\subset V_i$  and $W_n=\vide$. It follows that   $W_i$ is
finite, hence $\codim(W_i)>i$ for each $i>0$, so that the sheaf
$\alpha_*(E_f\otimes\omega_X)$ on $A$ is $M$-regular, hence \cgg. It is
therefore ample by Corollary \ref{coro}, and, since $\alpha$ is finite,
so are $\alpha^*\bigl(\alpha_*(E_f\otimes\omega_X)\bigr)$ and its
quotient  $E_f\otimes\omega_X$ (\S\ \ref{amp}).\end{proof}

In the following remarks, we keep the   hypotheses and notation of the theorem and its proof.

\begin{rema}  The
proof of the theorem shows that the sheaf $\alpha_*(E_f\otimes\omega_X)$  is 
\cgg.  In particular, if $f$ is not an isomorphism, $E_f\otimes\omega_X$ has nonzero
 sections, hence $p_g(Y)>p_g(X)$.
 \end{rema}

\begin{rema} The simplicity of the \av\ in the theorem is  essential: if $B$ is an \av\ and
 $g=(f,\Id_B):Y\times B\to X\times B$, we have $E_g=p^*E_f$, where  $p:X\times
B\to X$ is the first projection, hence $E_g\otimes\omega_{X\times B}=p^*(E_f \otimes\omega_X)$
 is not ample if $B$ is nonzero. The locus $W_i$ for $g$  contains
 $\PA\times
\{0\}$ for
$i\le\dim(B)$;  in particular, for $i=\dim(B)$, it is an abelian subvariety
 of codimension $i$ of $\Pic^0(A\times B)$.
 
\end{rema}

\begin{rema} If $X$ is not an \av, $ \omega_X$ is already ample (see, \eg, \cite{de}, Th\'eor\`eme 6.9)
 and one can show that the hypothesis that $f$ does not factor through a nontrivial connected abelian Galois \'etale
 covering of $X$ is unnecessary. If $X$ is a (simple) abelian variety, any finite cover $Y\to X$ factorizes as $Y\stackrel{f}{\to} X'\stackrel{\rho}{\to}X$ where $\rho$ is an isogeny and $f$ satisfies the hypotheses of the theorem.%
 \end{rema}

\begin{rema} Assume $X=A$ and let $d$ be the degree of $f$. For all  $i\ge
d-1$, the set $W_i$ is empty, \ie,
$$H^i(A,E_f\otimes P_\xi)=0\qquad{\rm for\ all\ \ }\xi\in\PA,
$$
by Le Potier's vanishing theorem
 (\cite{laz}, Theorem 7.3.5).
This does not hold in general for $0\le i< d-1$, as shown
by the following example.
Take
 an elliptic curve
 $C$, with origin $o_C$. Let $L$
 be a very ample
line bundle on $A$ and let  $Y\subset C\times A$ be a general  (smooth) element
of $|\cO_C((n+1)o_C) \boxtimes L |$. Following the proof of  \cite{laz}, Lemma
6.3.43, one sees that the second projection $f:Y\to A$ is finite (of degree $d=n+1$).
  By the Lefschetz theorem, the induced morphism
$$H^{n-i}(C\times A,\cO_{C\times A})\to H^{n-i}(Y,\cO_Y)
$$
is bijective for $i>0$ and injective for $i=0 $. In particular, $H^{n-i}(f,\cO)$ is
not surjective for $0\le i< n$, hence $0\in W_i$, \ie,
$$H^i(A,E_f)\ne 0\qquad{\rm for\ all\ \ }0\le i< d-1=n.
$$
In particular, $H^{n-1}(A,E_f)\ne 0 $, and it follows from \cite{muk},
Proposition
2.7, that  the $M$-regular vector bundle $E_f$ does not satisfy Mukai's
condition WIT${}_0$ when  $n>1$ (sheaves that satisfy condition  WIT${}_0$ are
$M$-regular).
\end{rema}

\backmatter

\end{document}